\documentclass[12pt]{article}
\usepackage{amsmath,amssymb,amsfonts,amsthm,url,xspace}
\newtheorem{theorem}{Theorem}

\newtheorem{lemma}[theorem]{Lemma}

\theoremstyle{definition}
\newcommand{\remark}{\medskip\noindent\textbf{Remark: }}
\def\endremark{\medskip}

\newcommand{\R}{\mathbb{R}}

\newcommand{\E}{\mathbb{E}}

\newcommand{\Var}{\operatorname{Var}}

\newcommand{\old}[1]{}

\newcommand{\tref}[1]{Theorem~\ref{thm:#1}}
\newcommand{\lref}[1]{Lemma~\ref{lem:#1}}

\newcommand{\cref}[1]{Corollary~\ref{cor:#1}}

\renewcommand{\th}{\ensuremath{^{\text{th}}}\xspace}
\def\rcs $#1: #2 ${\expandafter\def\csname rcs#1\endcsname {#2}}
\rcs $Date: 2004-10-11 15:21:29-07 $

\title{Balanced Boolean functions that can be evaluated \\ so that every input bit is unlikely to be read}
\author{Itai Benjamini\thanks{Weizmann Institute} \and Oded Schramm\thanks{Microsoft Research} \and David B. Wilson\footnotemark[2]}
\date{}
\begin{document}

\maketitle
\begin{abstract}
  A Boolean function of $n$ bits is balanced if it takes the value
  $1$ with probability $1/2$.  We exhibit a balanced Boolean function
  with a randomized evaluation procedure (with probability $0$ of
  making a mistake) so that on uniformly random inputs, no input bit
  is read with probability more than $\Theta(n^{-1/2} \sqrt{\log n})$.
  We give a balanced monotone Boolean function for which the
  corresponding probability is $\Theta(n^{-1/3} \log n)$.  We then
  show that for any randomized algorithm for evaluating a balanced
  Boolean function, when the input bits are uniformly random, there is
  some input bit that is read with probability at least
  $\Theta(n^{-1/2})$.  For balanced monotone Boolean functions, there
  is some input bit that is read with probability at least
  $\Theta(n^{-1/3})$.
\end{abstract}

\section{Results}

Suppose that a randomized algorithm evaluates a Boolean function
$f$ on $n$ input bits.  When the input bits are uniformly random, let
$\delta_i=\Pr[\text{bit $i$ gets read}]$, where the probability is
over the randomness of the input as well as the internal randomness of
the algorithm. 
Let $\delta=\max_i \delta_i$ be the maximum probability
that a particular bit is read.  How small can $\delta$ be for a
balanced Boolean function (one that takes on the values $0$ and $1$
equally often)?  What if the Boolean function is monotone
(i.e., $f(x)\ge f(y)$ whenever $x_i\ge y_i$ for all $i=1,2,\dots,n$)?

An obvious lower bound for $\delta_i$ is
$I_i(f)$, the influence of the bit $i$.
  Recall that the
influence $I_i=I_i(f)$ of the $i$\th input bit on $f$ is defined to be the
probability, for a uniformly random input, that changing the bit
changes the value of the function.  Many
readers are familiar with the majority function and the fact that the
influence of each bit is small, only $O(n^{-1/2})$.  But on most
inputs the numbers of $0$'s and $1$'s are nearly balanced, and any
algorithm that realiably evaluates majority will typically read
$\Theta(n)$ of the input bits, so that $\delta=\Theta(1)$.

There are a couple of other balanced Boolean functions that some
readers may think of.
The tribes function on $n=m\,2^m$ input bits is defined by partitioning the
set of bits into ``tribes'' of size $m$, and the function
takes the value $1$ if and only if there is at least one tribe where
all the bits are $1$.
  The influence
of each bit on the tribes function is only $\Theta(\log n/n)$ \cite{BenorLinial} (which is as small as $\max_i I_i(f)$ can be \cite{kkl}),
but any algorithm that
reliably evaluates the tribes function will typically read
$\Theta(n/\log n)$ of the input bits, so that $\delta=\Theta(1/\log
n)$.  For the dictatorship function, where one bit determines the
output, very few bits need to be read, but the dictator bit needs to be
read, so $\delta=1$.

The tribes function and the majority are two examples of symmetric functions;
they are invariant under a group acting transitively on the input bits.
Many other boolean functions of interest are symmetric.
Examples include recursive majority and
percolation crossings in a torus.
For symmetric functions
there is an algorithm $A$ computing $f$ that reads on average $m$ bits
if and only if there is an algorithm $A'$ computing $f$ with
$\delta\le m/n$. (To go from $A$ to $A'$, one just permutes the bits
by a uniformly-random element from the automorphism group before applying
$A$. The other direction is obvious, since the expected number
of bits that $A'$ reads is at most $\delta\,n$.)
Thus, for symmetric functions, estimates on $\delta$ are equivalent
to estimates on the expected number of bits that need to be read.

In the
next section we give some nearly optimal constructions of Boolean
functions that may be evaluated so that every input bit is unlikely to
be read:
\begin{theorem}
\label{thm:upper}
\begin{enumerate}
\item There is a balanced Boolean function with an algorithm that
  always correctly evaluates it on any input and for which
  $\delta=\Theta(n^{-1/2} \sqrt{\log n})$.
\item There is a balanced Boolean function with an algorithm that
  correctly evaluates it on most inputs most of the time and for which
  $\delta=\Theta(n^{-1/2})$.
\item There is a balanced monotone Boolean function with an algorithm
  that always correctly evaluates it on any input and for which
  $\delta=\Theta(n^{-1/3} \log n)$.
\item There is a balanced monotone Boolean function with an algorithm
  that correctly evaluates it on most inputs most of the time and for
  which $\delta=\Theta(n^{-1/3})$.
\end{enumerate}
\end{theorem}

The constructions of \tref{upper} are optimal except possibly for the
factors of $\sqrt{\log n}$ and $\log n$, as the lower bounds that we
prove in Sections~\ref{sec:lowbd} and~\ref{sec:mtnlowbd} show:
\begin{theorem}
\label{thm:lower}
\begin{enumerate}
\item If an algorithm correctly computes a balanced Boolean function
  $f:\{-1,1\}^n\to\{-1,1\}$
  on most inputs most of the time, then $\delta\geq\Theta(n^{-1/2})$.
More generally and more precisely,
 $\Pr[\text{algorithm is wrong}]\ge \frac18\Var(f)-\frac14n\delta^2$,
 holds regardless if $f$ is balanced.
\item If an algorithm correctly computes a balanced monotone Boolean
  function on all inputs all of the time, then
  $\delta\geq\Theta(n^{-1/3})$.
  More generally and more precisely,
  if $f$ is monotone but not necessarily balanced,
  then $\Var(f)\le \delta^{3/2}\,n^{1/2}$.
\end{enumerate}
\end{theorem}

\section{Constructions}
The constructions are based on directed percolation on certain graphs,
and were inspired in part by coupling-from-the-past on generic Markov
chains \cite{propp-wilson:unknown-markov-tree} and in part by Radford
Neal's circular coupling \cite{neal}, though background on these
topics is not required to understand the constructions.

\subsection{Directed percolation on the wrapped extended butterfly}

The examples we construct are based on a directed graph.
There are various different choices that might work here, and the
framework we have chosen, the wrapped extended butterfly, offers
an explicit description and reasonably clear proofs.

Define the wrapped extended butterfly $\Omega_{H,W}$ to be a directed
graph with $H\times W$ vertices, where $H$ is even and often a power
of $2$.  The coordinates of a vertex are $(h,t)$ where $h$ and $t$ are
integers with $0\leq h<H$ and $0\leq t<W$.  Each vertex has in-degree
two and out-degree two.  The vertex $(h,t)$ has directed edges leading
to $(2 h \bmod H,t+1\bmod W)$ and $(2h+1\bmod H,t+1\bmod W)$.  When
$H=2^d$ and $W=d$, the undirected unwrapped version of this graph is
sometimes known as the ``omega network'' or ``shuffle network''
\cite[Chapter 3.8.1]{MR92j:68041}, which is isomorphic to the usual
butterfly \cite[Chapter 3.2.1]{MR92j:68041} that is sometimes used in
parallel computing architectures.  The extended butterfly (when
$H=2^d$ and $W= 5 d$) has also been used in the construction of
holographic proof systems \cite{195132}.  We will for convenience take
$H=2^d$, but for our constructions $W$ will be much larger,
$W=\Theta(d 2^d)$ for one construction and $W=\Theta(2^{d/2})$ for the
other.  There are many graphs that we could have used in our
constructions, but with the wrapped extended butterfly, the
calculations come out nicely.

We will refer to the set of vertices with last coordinate $t$ as the
``$t$\th time slice'' (where ``time'' is periodic), and the vertices
with the same first coordinate as ``points''.

We will consider two ensembles of random subgraphs of the wrapped
extended butterfly.  In the non-monotone ensemble, which we use for the
non-monotone construction, each vertex has exactly one of its two
out-going edges included in the subgraph, and an input bit determines
the least significant bit (i.e., the parity)
of the first coordinate of the destination vertex.  We think of the $n=H
W$ input bits as being associated with the vertices.  In the monotone
ensemble, which we use for the monotone construction, both possible
out-going edges from a vertex have an input bit to determine whether
or not the edge is present in the subgraph.  Here the expected
out-degree of each vertex is $1$.  A total of $n=2 H W$ input bits
(associated with the edges) are used.

In the various constructions we will be interested in the directed
cycles that may exist in this random subgraph defined by the input
bits.  A Las Vegas algorithm may determine the cycles that exist in
the random subgraph by picking a uniformly random time slice and
following all the paths forward in time until the starting time slice
is reached.  No additional bits of the input need to be read.  Many of
the paths merge early on, so following them all does not require as
many reads as it might at first appear.  As we shall see later, a
Monte Carlo algorithm that works
most of the time for most inputs may instead follow
the paths starting from a smaller random set of vertices.

\subsection{The nonmonotone random subgraph ensemble}

Since the subgraphs in the nonmonotone ensemble have out-degree $1$ at
each vertex, there must be at least one directed cycle. 
The Boolean function can be any sort of
balanced function of these cycles.  For example, we could take the
lexicographically smallest vertex that is in a cycle and in time slice
$0$, and take its bit.
\begin{lemma}
  Provided $H=2^d$ and $W> d$, this Boolean function is exactly balanced.
\end{lemma}
\begin{proof}
  Because the wrapped extended butterfly is symmetric under the
  operation that flips the interpretation of the bits at a given time
  slice, a random subgraph is isomorphic to the subgraph obtained by
  flipping the bits at a time slice and suitably permuting the bits at
  the next $d$ time slices.
\end{proof}
\remark If we wish the Boolean function to be symmetric in its input
bits while remaining exactly balanced, we can do this via a trick that
places four input bits $b_{h,t,0},\ldots,b_{h,t,3}$ at each vertex
$(h,t)$.  The parity ($b_{h,t,0}\oplus\cdots\oplus b_{h,t,3}$) of
these four bits determines which outgoing edge from vertex $(h,t)$ is
used.  These four bits also determine another bit $b'_{h,t}$
that is symmetric in
the four bits and independent of parity: this bit $b'_{h,t}$ is $1$
iff there are exactly one or two ones among
$(b_{h,t,0},b_{h,t,1},b_{h,t,2},b_{h,t,3})$ and these are consecutive
in the cyclic order.
We can take our
Boolean function $f$ to be the XOR of these $b'_{h,t}$
bits at all vertices that lie
in a cycle.  For any $H=2^d$ and $W>d$, the Boolean function is
symmetric in the input bits and exactly balanced.
\endremark

\begin{lemma}\label{lem:nonmonbd}
  In the nonmonotone ensemble, the probability that the Las Vegas
  algorithm reads a bit is at most $O(W^{-1} \log H + H^{-1})$.
\end{lemma}
\begin{proof}
  Let $p_t$ be the number of vertices whose bits are read
  in a time slice $t$ time units after the
  initial time slice.  If we run the $p$ particles independently
  forward $\log_2 H$
  time steps from these locations,
  then their resulting locations would be independent
  and uniformly random. Under this independent random
  walks dynamics, when two particles
  reach a vertex, their outgoing edges may be different.
  If instead we consider coalescing dynamics, each vertex has one
  outgoing edge (and these outgoing edges are chosen uniformly independently),
  and all the particles arriving at the vertex always use that same
  edge. We may couple the coalescing dynamics with the independent random walks
  dynamics in such a way that the occupied locations at the end of
  the walks in the coalescing dynamics is a subset of the locations
  in the independent random walks dynamics.
  Suppose that two
  independent uniformly random particles in the same
  time slice do non-coalescing independent
  random walks on the wrapped extended butterfly for $s$ steps each,
  where $s<W$, and let $N$ be
  the number of collisions, that is, vertices visited by both.
  Then  $\E[N] = s/H$. Moveover, given that the particles visited
  the same vertex, the probability that they will be at the same
  vertex $t$ time slices later is $\max\{2^{-t},1/H\}$.
  Consequently, $\E[N|N>0] < 2+s/H$, and therefore
  $\Pr[N>0]=\E[N]/\E[N|N>0] > s/(2H+s)$. 
  Regardless of what the two particles do, the probability
  that a third particle (started from an independent and uniformly
  random location in their initial time slice)
  collides with at least one of them is at most $2
  s/H$.  Thus the probability that a given pair of particles collide
  but that neither one collides with another particle is at least
  $$ \frac{s}{2H+s} \left[1-\frac{2s}{H}\right]^{p_t}.$$
  Under the
  coalescing dynamics, the probability of this event can only be
  larger.  The number of such events is at most the reduction in
  the number of particles, so we find
  $$
  \E[p_{t+\log_2 H + s}|p_t] \leq p_t -
  \binom{p_t}{2} \frac{s}{2H+s} \left[1-\frac{2s}{H}\right]^{p_t}.$$
  Taking $s\approx H/p_t$, we see that the expected reduction
  (that is, 
  $\E[p_t- p_{t+\log_2 H + s}|p_t ]$) is
  $\Theta(p_t)$.  Let $T_k$ be the number of time slices for which $2^k
  < p_t \leq 2^{k+1}$.  From the above, we see that
  $\E[T_k]=O(H/2^k+\log H)$.  Thus $\sum_k 2^{k+1} \E[T_k] \leq O(H
  \log H)$.  After the paths have all coalesced, there is one path that
  continues for $O(W)$ steps.  Since the expected number of input bits
  read is $O(H\log H + W)$, and each of the $H W$ input bits is read
  with the same probability, we obtain the desired bound on $\delta$.
\end{proof}

The optimal choice is $H =\Theta(\sqrt{n/\log n})$ and
$W=\Theta(\sqrt{n\log n})$, giving $\delta = O(n^{-1/2} \sqrt{\log n})$,
which gives part 1 of Theorem 1.

\bigskip
To obtain part 2 of Theorem 1 we instead take $W=H$.  The Las Vegas
algorithm would have $\delta=O(n^{-1/2} \log n)$, but we can save the
factor of $\log n$ by settling for an algorithm that reports the
correct answer most of the time on most inputs.  Rather than locating
the set of cycles by trying all possible starting points at a random
time slice and following them forward $W$ steps, we instead pick only
$m$  points from the uniformly selected time slice $t_0$, and
follow them forward until they cycle.
Here, $m$ is large, but depends only on the required reliability of
the algorithm. Given $t_0$, the $m$ points are selected at random,
uniformly and independently.

\begin{lemma}\label{lem:nr}
  Assuming $W=H$, the probability that this algorithm reads any
  specific bit is $O(m/H)$
\end{lemma}
\begin{proof}
We prove this in the case $m=1$, which is clearly sufficient.
Moreover, by symmetry, it is enough to show that the expected
number of vertices visited is $O(W)$.
We start with a single particle, and follow its path.
Let $X$ be the number of vertices visited.
If the path did not enter its previously visited vertices after $s_0$ steps,
where $s_0\ge W$, then the conditioned probability  that it will
hit itself in the next $s+\log_2 H$ steps is at least
$s/(2H+s)$, by the argument giving the lower bound for $\Pr[N>0]$
in the proof of \lref{nonmonbd}. Taking $s=H-\log_2 H$,
say, we find that $\Pr[X>(j+1)W|X>jW]<3/4$ for every $j=1,2,\dots$.
Consequently, $\E[X]=O(W)$.
\end{proof}

To complete the argument, we need to show that with probability tending
to $1$ as $m\to\infty$ (uniformly in $H$), the algorithm finds all
open cycles. This follows from the following lemma.

\begin{lemma}
  Assuming $W=H$, the expected number of open cycles that are
  undiscovered by the algorithm is $O(1/(m+1))$.
\end{lemma}

\begin{proof}
Let $\{v_0,v_1,\ldots,v_m\}$
be a set of $m+1$ uniformly chosen vertices in the $t_0$
time slice, which are chosen independently given $t_0$.
Let $\Gamma$ be the union of the paths of length
$W/2$ starting at $\{v_1,v_2,\ldots,v_m\}$.
We are going bound the probability for the event $\mathcal A$ that
$v_0$ is on an open cycle that does not intersect $\Gamma$.
Let $\gamma$ be the path of length $W/2$ starting
at $v_0$.

If we were to follow all $H$ particles at time slice $t_0$ forward in
time, let $\tau_j$ be the number of time steps before there are at
most $j$ particles.  Because $\Pr[\tau_j\geq\lceil
2\E[\tau_j]\rceil]\leq 1/2$ and successive blocks of $\lceil
2\E[\tau_j]\rceil$ time slices are independent, $\Pr[\tau_j\geq\ell\lceil
2\E[\tau_j]\rceil]\leq 1/2^\ell$.  Provided $j$ is not too large
($j\leq H/(\log H \log\log H)$), our bound for $\E[T_k]$ in \lref{nonmonbd}
implies that $\E[\tau_j]=O(H/j)$, so we can take $\ell=\Theta(j)$ to find
that $\Pr[\tau_j>H/2]$ decays geometrically in $j$.

Let $Y$ be the number of particles starting
at $\{v_0,\ldots,v_m\}$ that have not merged with any
other particle after $W/2$ steps;
we have $\E[Y]\le O(1)$.
By symmetry
 $\Pr[\gamma\cap\Gamma=\emptyset|Y]= Y/(m+1)$, so in fact
$\Pr[\gamma\cap\Gamma=\emptyset]\le O(1/(m+1))$.
Given $\gamma$ and $\Gamma$, the probability
that the continuation of $\gamma$ hits $v_0$ on
first return to the $t_0$ time slice is 
exactly $1/H$.
Moreover, the same argument as in 
\lref{nr} shows that given $\gamma$ and $\Gamma$
the probability that the continuation of $\gamma$ hits
$v_0$ at the $j$\th return to time slice $t_0$ (but not sooner) is
bounded by $O((3/4)^{j}/H)$.
Thus, we have
$\Pr[\mathcal A]\le O(1/((m+1)H))$.
Since the number of open cycles that are disjoint from $\Gamma$
is at most the number of vertices in the $t_0$ time slice that
are on such cycles, this proves the lemma.
\end{proof}

\subsection{The monotone random subgraph ensemble}

Since the subgraphs in the monotone random subgraph ensemble have
random out-degree, we are not assured that there is a directed cycle,
and when $W\gg\sqrt{H}$ it is even unlikely for there to be a directed
cycle.  This suggests that we take our Boolean function to be the
existence of a directed cycle, or what turns out to be a little
simpler to analyze, the existence of a directed cycle of length $W$.

\begin{lemma}
  In the monotone ensemble, the probability that the Las Vegas reads a
  bit is at most $(4+o(1)) W^{-1} \log W$.
\end{lemma}
\begin{proof}
  Let $v$ be the starting vertex of the edge assigned to the given
  input bit, and suppose that the initial time slice was $t$ steps
  before $v$. Let $G_v$ be the subgraph of the wrapped extended butterfly
  which consists of all directed paths of length at most $W$ terminating at $v$.
  Each edge of $G_v$
  occurs in the random subgraph with probability $1/2$.
  There is an obvious graph covering of $G_v$ by the binary tree of depth $W$.
  An upper bound for the probability that $v$ is explored conditioned on
  $t$ is therefore the probability that the root of the binary tree is
  connected to some vertex at level $t$ in a subgraph where each
  edge is included with probability $1/2$ independently.
  (See, e.g.,~\cite{campanino-russo} or~\cite[Theorem 1]{bs-pyond}.)
    Let $S_t$ be this latter probability.  We
  have $S_{t+1}=S_t-\frac14 S_t^2$, so for large $t$, $S_t\sim 4/t$.
  The probability that a vertex is explored is upper bounded by
  \begin{align*}\frac1 W\sum_{t=0}^{W-1} S_t &= (4+o(1)) \frac{\log W}{W}. \qedhere\end{align*}
\end{proof}

In the monotone ensemble, let us consider the cycles that wind around
exactly once (cycles with length $W$).  Assuming $W\geq d$, these
cycles are in bijective correspondence with bit strings of
length $W$: given a cycle, the bit string is formed by the parities of
the vertical coordinates, and given a bit string, the vertical
coordinates are determined by groups of $d$ consecutive bits (in circular
order).  Let
$N_{H,W}$ be the number of such cycles that occur in the graph.  Since
there are $2^W$ possible such cycles, each occuring with probablity
$2^{-W}$, we have $\E[N_{H,W}]=1$.

\remark
As we shall see, the interesting regime is when $W=\Theta(\sqrt{H})$.
It is not hard to show
$\lim_{1\ll W\ll \sqrt{H}}\Pr[N_{H,W}>0]= 1-1/e$ and
$\lim_{1\ll \sqrt{H} \ll W}\Pr[N_{H,W}>0]=0$.  Presumably
$\lim_{H\to\infty}\Pr[N_{H,c\sqrt{H}}>0]$ is a
continuous function of $c$, in which case there is some particular
positive value of $c$ that we can pick so that
$\Pr[N_{H,c\sqrt{H}}>0]\doteq 1/2$,
so that our Boolean function is nearly balanced. 
$\lim_{H\to\infty}\Pr[N_{H,c\sqrt{H}}>0]$ is an
interesting function of $c$, but not one that we shall explore here.
For our purposes it suffices to use a second moment estimate.
\endremark

To compute $\E[N_{H,W}^2]$, we are interested in pairs of circular
bit strings.  Some pairs of cycles share edges, and for these pairs
the probability that both cycles occur will be larger than $2^{-2 W}$.
Define a merge time to be a time slice where the two cycles
coincide on the outgoing edge but not the incoming edge, and a split
time to be a time slice where the two cycles coincide on the incoming
edge but not the outgoing edge.  The splits and merges alternate, and
after a split, the next merge cannot occur within the next $d$ time
slices, but there are no other constraints on when the splits and
merges may occur.  Suppose there are $\ell$ splits and $\ell$ merges
--- these occur at distinct times, so there are at most
$2\binom{W}{2\ell}$ ways to select the times during which the two
cycles share an edge.
(When $\ell=0$, there is one way to select the merge and split times,
but still two ways to decide which edges agree between the cycles.)
Given the
times during which the cycles share edges, how many pairs of cycles
satisfy these constraints?
There are $2^W$ ways to select the first cycle, each occuring with probability
$2^{-W}$.  Every merge specifies
the preceding $d+1$ bits in the second cycle, and the bit for every shared
edge as well as every edge at a split time is also specified. 
The number of
ways to pick the second string is at most $2^{W-S}/(2 H)^\ell$,
where $S$ is the number of shared edges between the two cycles.
The probability that the
second cycle occurs given that the first one does is $2^{-W+S}$.  The
expected number of pairs of cycles with $\ell$ splits and $\ell$
merges is at most
$$ 2 \frac{W^{2\ell}}{(2\ell)!}\frac{1}{(2 H)^\ell}.$$
We may sum over all $\ell$ to conclude
$
\E[N_{H,W}^2] \leq \exp(W/\sqrt{2H}) +\exp(-W/\sqrt{2H})$.
Since
$\Pr[N_{H,W}>0]\geq\E[N_{H,W}]^2/\E[N_{H,W}^2]$, we have
$$\Pr[N_{H,W}>0]\geq \frac{1}{\exp(W/\sqrt{2H}) +\exp(-W/\sqrt{2H})}.$$

We can pick that value of $c$ for which $1/(e^c+e^{-c})=1-1/\sqrt{2}$
($c\doteq 1.12838$) and set $W=\lfloor c\sqrt{2 H}\rfloor$.
Instead of asking about the existence of a cycle
of length $W$, we instead ask if there is a cycle that passes through
a suitable set of vertices at the first time slice.  The probability
that a particular vertex is part of a cycle is at most $1/H$, the
probability that a vertex is part of a cycle while the
lexicographically smaller ones aren't will be even smaller, so by
adjusting the size of the set of vertices we can ensure $1-1/\sqrt{2} \leq
\Pr[\text{suitable cycle exists}] \leq 1-1/\sqrt{2} + 1/H$.  If a
cycle passes through the last node in the suitable set, then call it a
marginally suitable cycle, otherwise call it completely suitable:
$\Pr[\text{completely suitable cycle exists}]\leq 1-1/\sqrt{2}$.  We
may repeat this experiment twice
(with new bits for the second experiment)
while taking $H=\Theta(n^{2/3})$ and
$W=\lfloor c \sqrt{2 H}\rfloor$ in both experiments.  If there is a
completely suitable cycle in either experiment, the Boolean function
takes the value $1$, if there is no suitable cycle in either
experiment, the Boolean function takes the value $0$.  In the third
scenario, where there is a marginally suitable cycle but not a completely
suitable cycle, we may tweak the definition of the Boolean function so
that it is monotone and exactly balanced.  This third scenario occurs
with probability $O(n^{-2/3})$, so it can increase $\delta$ by at most
$O(n^{-2/3})$, so we still have $\delta\leq O(n^{-1/3} \log n)$,
giving us part 3 of Theorem 1.

\bigskip
For part 4 of Theorem 1, we can use the same Boolean function that we
used for part 3, and use the same choice of
$W$ and $H$.  As before, we can save the factor of $\log n$ by
choosing a typically smaller collection of starting points.
Let $m$ be a large constant that depends only on the required reliability
of the algorithm.
The algorithm chooses a
random  set $S$ of vertices in the extended butterfly,
where each vertex is in $S$ with probability $m/H$, independently.
The algorithm then explores all the open paths of length at most
$2W$ starting at the vertices in $S$.

As above, $\delta$
can be bounded by a corresponding process on the binary tree of depth $2\,W$.
In this case, $\delta$ is bounded by the expected number of
\lq\lq selected\rq\rq\ vertices in the percolation component of the root,
where each vertex is selected with probability $m/H$ independently
from other vertices and from the percolation process.
Since the expected number of vertices
in the percolation cluster of the root and at distance $s$ from the root is $1$
(when $s\le 2\,W$), we get $\delta\le 2\,W\,m/H=O(m\,n^{-1/3})$.

The proof of the theorem is completed by showing that
this algorithm finds all open cycles of length $W$ with high probability.
This follows from the following lemma.

\begin{lemma}
  Assuming $\sqrt m<W=\Theta(\sqrt H)$, the expected number of open
  cycles of length $W$ unexplored by this algorithm is at most
  $c^{-\sqrt m}$ for some constant $c\in(0,1)$.
\end{lemma}

\begin{proof}
Fix some cycle $\gamma$ of length $W$. Let $\mathcal A$ be the event
that $\gamma$ is open and undetected by the algorithm.
Set $W_0=\lfloor W/\sqrt m\rfloor$, and
let $\gamma_0$ be a subpath of $\gamma$ with $W_0$ vertices.
Let $N$ be the number of open paths of length at most
$W_0$ that start at a vertex in $S$, end at a vertex in
$\gamma_0$, and are otherwise disjoint from $\gamma$.
We will now prove that $\Pr[N>0]=\Theta(1)$ by a second
moment argument.

Fix some $k\in\{0,1,2,\dots,W_0\}$.
The number of paths of length $k$ that end at a vertex in
$\gamma_0$ but are otherwise disjoint from $\gamma$ is $\Theta(W_0\,2^k)$.
Thus 
$$
\E[N]=\sum_{k=0}^{W_0}
\Theta(W_0\,2^k)\,2^{-k}\,\frac{m}{H}
=\Theta(1)\,.
$$
The estimation of the second moment is similar to the one done
above for $N_{H,W}$.  We can bound $\E[N^2]$ by considering separately
the four cases where the pair of paths have the same or different
starting points and the same or different ending points.  We enumerate
the pairs of paths by tracing them backwards from $\gamma_0$.  Let
$\ell$ be the number of splits (merges in reverse), and let $k_1$ and
$k_2$ be the lengths of the two paths.  We find $\E[N^2]$ is at most
\begin{multline*}
  \sum_{\smash{k_1=0}}^{W_0} W_0 2^{k_1} 2^{-k_1} \frac{m}{H}\times
  \Bigg[
  \overbrace{\sum_{k_2=0}^{W_0} W_0 2^{k_2} 2^{-k_2} \frac{m}{H} 
  \sum_{\ell\ge0}\binom{W_0}{2\ell} \frac{1}{H^\ell}}
  ^{\text{different start, different end}} +
  \overbrace{\sum_{k_2=0}^{W_0} 2^{k_2} 2^{-k_2} \frac{m}{H} 
  \sum_{\ell\ge0}\binom{W_0}{2\ell+1} \frac{1}{H^\ell}}
  ^{\text{different start, same end}} +\\
  \overbrace{\sum_{k_2=0}^{W_0} 2^{k_2} 2^{-k_2}
  \sum_{\ell\ge1}\binom{W_0}{2\ell-1} \frac{1}{H^\ell}}
  ^{\text{same start, different end}} +
  \overbrace{\sum_{k_2=k_1} 2^{k_2} 2^{-k_2}
  \sum_{\ell\ge0}\binom{W_0}{2\ell} \frac{1}{H^\ell}}
  ^{\text{same start, same end}}
\Bigg] = O(1),
\end{multline*}
and thus $\Pr[N>0]=\Theta(1)$.

We may conclude that for every $k=0,1,\ldots,\lfloor\sqrt m/2\rfloor-1$
the probability that $\gamma$ is hit in a time slice
$t\in [2\,k\,W_0,(2\,k+1)\,W_0)$ by an open path of length
at most $W_0$ starting in $S$ is $\Theta(1)$.
Since these events are independent, the probability that
$\gamma$ is not visited by an open path of length
at most $W_0$ starting in $S$ is $c^{-\sqrt m}$ for some constant
$c\in(0,1)$.  The proof is
complete by noting that there are $2^W$ possible cycles $\gamma$, and
each is open with probability $2^{-W}$.
\end{proof}

\section{Lower bound for balanced Boolean functions}
\label{sec:lowbd}

Suppose that a randomized algorithm approximately computes a Boolean
function.  Let $A(r,z)$ denote the output of the algorithm when it
uses coins $r=r_1 r_2\cdots$ on input $z=z_1z_2\cdots z_n$.  As usual
we let $\delta_i=\Pr[\text{bit $i$ gets read}]$.  Consider two
independent runs of the algorithm (i.e., using independent coins
$r=r_1 r_2\cdots$ and $s=s_1 s_2\cdots$) on independent inputs
$x=x_1x_2\cdots x_n$ and $y=y_1y_2\cdots y_n$.  Let $N$ be the number
of bit positions that are read by both of these independent runs:
$$\Pr[N>0]\leq \E[N] = \sum_i \delta_i^2 \leq n \delta^2.$$
We define $z=z_1z_2\cdots z_n$ by
$$z_i=\begin{cases}x_i & \text{if algorithm with coins $r$ and input
    $x$ reads bit $i$}\\ y_i & \text{otherwise.}\end{cases}$$
The
vector $z$ is uniformly random and independent of $r$ and $s$.  Of
course $A(r,z)=A(r,x)$.  In the event that $N=0$, we have
$A(s,z)=A(s,y)$.  Since $A(r,x)$ and $A(s,y)$ are i.i.d.,
\begin{multline*}
\Pr[A(r,z)=A(s,z)]
\le \Pr[A(r,x)=A(s,y)] + \Pr[N>0] \\
 \leq \Pr[A(r,x)=-1]^2 + \Pr[A(r,x)=1]^2 +
n\delta^2.
\end{multline*}
Set $w= \Pr[A(r,x)\ne f(x)]$, $p=\Pr[f(x)=1]$
and $p'=\Pr[A(r,x)=1]$.
If $A(r,z)\ne A(s,z)$, then either $A(r,z)\ne f(z)$ or $A(s,z)\ne f(z)$.
Thus,
$$
w
\ge \frac12\bigl( 1-{p'}^2 - (1-p')^2 - n\delta^2\bigr)
=
p'-{p'}^2 -\frac 12 n\delta^2\,.
$$
Now note that the absolute value of the $p$-derivative of $\Var(f)/4=p-p^2$
is bounded by $1$. Thus, $p'-{p'}^2\ge\frac14\Var(f)-|p-p'|\ge \frac14\Var(f)-w$.
Consequently,
$$
2w \ge \frac 14\Var(f)-\frac 12n\delta^2\,,
$$
as required.

\section{Lower bound for balanced monotone Boolean functions}
\label{sec:mtnlowbd}

Suppose that $f:\{-1,1\}^n\to\R$ is some function, and
$A$ is a randomized algorithm calculating $f$ (exactly and always).
Recall the definition of the Fourier coefficients
$$
\hat f(S):= \E\Bigl[ f(x)\prod_{i\in S} x_i\Bigr]\,,
\qquad\qquad S\subset\{1,2,\dots,n\}\,.
$$
When $f$ is monotone and takes values in $\{-1,1\}$, we clearly have
$$
I_i(f)=\hat f(\{i\})\,.
$$
We will need the inequality
\begin{equation}\label{eqn:i}
\sum_{i=1}^n\hat f(\{i\}) \le \sqrt{n\,\delta}\,.
\end{equation}
This inequality holds for $f$ taking values in $\{-1,1\}$,
even if $f$ is not monotone.
It is obtained by combining the first two displayed inequalities
in the proof of Theorem 1 in~\cite{odonnell-servedio}
(although the algorithms discussed there are deterministic, the proof
applies to random algorithms as well).
Alternatively, the case $k=1$ of the inequality
\begin{equation}\label{eqn:ss}
\sum_{|S|=k}\hat f(S)^2\le \delta\,k\,\|f\|_2^2
\end{equation}
from~\cite{schramm-steif} gives
$$
\Bigl(\sum_{i=1}^n \hat f(\{i\})\Bigr)^2
\le n\,\sum_{i=1}^n \hat f(\{i\})^2 \le n\, \delta\,,
$$
implying~\eqref{eqn:i}.
(Although~\eqref{eqn:ss} does not assume that $f$ is monotone or
boolean, in the last step we assumed that $f$ takes values
in $\{-1,1\}$ to drop the factor $\|f\|_2^2$.)

Another inequality that we need to quote (valid for boolean but
not necessarily monotone $f$) is
\begin{equation}\label{eqn:osss}
\Var(f)\le \sum_{i=1}^n \delta_i\,I_i(f)\,.
\end{equation}
See~\cite{odonnell-saks-schramm-servedio}.

Now assume that $f:\{-1,1\}^n\to\{-1,1\}$ is monotone.
Since $I_i(f)=\hat f(\{i\})$,
inequalities~\eqref{eqn:osss} and~\eqref{eqn:i} give
$$
\Var(f)\le \delta^{3/2}\,n^{1/2},
$$
which proves part 2 of \tref{lower}.

\section*{Acknowledgements}
We thank Ryan O'Donnell, Mike Saks and
Rocco Servedio for useful discussions and for permission to use
inequality~\eqref{eqn:osss}.

\bibliographystyle{halpha}
\bibliography{ptr}

\begin{thebibliography}{OSSS04}

\bibitem[BOL89]{BenorLinial}
Michael Ben-Or and Nathan Linial.
\newblock Collective coin flipping.
\newblock In S.~Micali, editor, {\em Randomness and Computation}, pages
  91--115, New York, 1989. Academic Press.

\bibitem[BS96]{bs-pyond}
Itai Benjamini and Oded Schramm.
\newblock Percolation beyond {$\Z^d$}, many questions and a few answers.
\newblock {\em Electron. Comm. Probab.}, 1(8):71--82 (electronic), 1996,
  \url{http://www.math.washington.edu/~ejpecp/EcpVol1/paper8.abs.html}.

\bibitem[CR85]{campanino-russo}
M.~Campanino and L.~Russo.
\newblock An upper bound on the critical percolation probability for the
  three-dimensional cubic lattice.
\newblock {\em Ann. Probab.}, 13(2):478--491, 1985.

\bibitem[KKL88]{kkl}
Jeff Kahn, Gil Kalai, and Nathan Linial.
\newblock The influence of variables on boolean functions (extended abstract).
\newblock In {\em 29th Annual Symposium on Foundations of Computer Science},
  pages 68--80, 1988.

\bibitem[Lei92]{MR92j:68041}
F.~Thomson Leighton.
\newblock {\em Introduction to Parallel Algorithms and Architectures: Arrays,
  Trees, Hypercubes}.
\newblock Morgan Kaufmann, San Mateo, CA, 1992.

\bibitem[Nea02]{neal}
Radford~M. Neal.
\newblock Circularly-coupled {Markov} chain sampling.
\newblock Technical Report 9910 (revised), Dept. of Statistics, University of
  Toronto, 2002.

\bibitem[OS04]{odonnell-servedio}
Ryan O'Donnell and Rocco Servedio.
\newblock On decision trees, influences, and learning monotone decision trees.
\newblock Technical Report CUCS-023-04, Columbia University, Dept. of Computer
  Science, 2004, \url{http://www1.cs.columbia.edu/~library/2004.html}.

\bibitem[OSSS04]{odonnell-saks-schramm-servedio}
Ryan O'Donnell, Mike Saks, Oded Schramm, and Rocco Servedio, 2004.
\newblock Manuscript.

\bibitem[PS94]{195132}
Alexander Polishchuk and Daniel~A. Spielman.
\newblock Nearly-linear size holographic proofs.
\newblock In {\em Proceedings of the twenty-sixth annual ACM symposium on
  Theory of computing}, pages 194--203. ACM Press, 1994.

\bibitem[PW98]{propp-wilson:unknown-markov-tree}
James~G. Propp and David~B. Wilson.
\newblock How to get a perfectly random sample from a generic {M}arkov chain
  and generate a random spanning tree of a directed graph.
\newblock {\em Journal of Algorithms}, 27:170--217, 1998.

\bibitem[SS04]{schramm-steif}
Oded Schramm and Jeff Steif.
\newblock Quantitative noise sensitivity and exceptional times for percolation,
  2004.
\newblock In preparation.

\end{thebibliography}

\end{document}